# ADAPTIVE BAYESIAN ESTIMATION USING A GAUSSIAN RANDOM FIELD WITH INVERSE GAMMA BANDWIDTH

By A. W. van der Vaart and J. H. van Zanten[1]

*Vrije Universiteit Amsterdam*

We consider nonparametric Bayesian estimation inference using a rescaled smooth Gaussian field as a prior for a multidimensional function. The rescaling is achieved using a Gamma variable and the procedure can be viewed as choosing an inverse Gamma bandwidth. The procedure is studied from a frequentist perspective in three statistical settings involving replicated observations (density estimation, regression and classification). We prove that the resulting posterior distribution shrinks to the distribution that generates the data at a speed which is minimax-optimal up to a logarithmic factor, whatever the regularity level of the data-generating distribution. Thus the hierachical Bayesian procedure, with a fixed prior, is shown to be fully adaptive.

**1. Introduction.** The quality of nonparametric estimators of densities or regression functions is well known to depend on the regularity of the true density or regression function. Given $n$ independent observations on a function of $d$ arguments that is only known to be $\alpha$-smooth, the precision of estimation is of the order $n^{-\alpha/(2\alpha+d)}$. Initially this was shown using estimators that depend explicitly on the regularity level $\alpha$, but later it was shown that the optimal rate can be achieved for all levels of regularity simultaneously. Estimators that are rate optimal for every regularity level are called *adaptive*. Cross validation, thresholding, penalization and blocking are typical methods to construct such estimators (see, e.g., [1, 2, 6, 10, 11, 12, 13, 19, 33, 34, 35, 37] and [42]).

Adaptive methods often employ a scale of estimators indexed by a bandwidth parameter and adapt by making a data-dependent choice of the bandwidth. Within a Bayesian context it is natural to put a prior on such a

Received July 2008; revised December 2008.
[1]Supported in part by the Netherlands Organization for Scientific Research NWO.

bandwidth parameter and let the bandwidth be chosen through the posterior distribution. In this paper we discuss a particularly attractive Bayesian scheme, and show that this yields estimators that are adaptive up to a logarithmic factor.

Our scheme employs a fixed prior distribution, constructed by rescaling a smooth Gaussian random field. There is some (but not much) freedom in the choice of Gaussian field and scaling factor. One possible choice is the squared exponential process combined with an inverse Gamma bandwidth. The *squared exponential process* is the centered Gaussian process $W = \{W_t : t \in \mathbb{R}^d\}$ with covariance function, for $\|\cdot\|$ the Euclidean norm on $\mathbb{R}^d$,

$$(1.1) \qquad EW_sW_t = \exp(-\|t-s\|^2).$$

The Gaussian field $W$ is well known to have a version with infinitely smooth sample paths $t \mapsto W_t$. To make it suitable as a prior for $\alpha$-smooth functions we rescale the sample paths by an independent random variable $A$ distributed as the $d$th root of a Gamma variable. As a prior distribution for a function on the domain $[0,1]^d$ we consider the law of the process

$$\{W_{At} : t \in [0,1]^d\}.$$

The inverse $1/A$ of the variable $A$ can be viewed as a bandwidth parameter. For large $A$ the prior sample path $t \mapsto W_{At}$ is obtained by shrinking the long sample path $t \mapsto W_t$ indexed by $t \in [0,A]^d$ to the unit cube $[0,1]^d$. Thus it employs "more randomness" and becomes suitable as a prior model for less regular functions if $A$ is large.

The effect of scaling the prior was already noted in [47], who showed (for $d=1$) that a deterministic scaling by the "usual" bandwidth $1/A = n^{-1/(2\alpha+1)}$ produces priors that are suitable models for $\alpha$-regular functions. The main contribution of the present paper is to show that a single inverse Gamma bandwidth gives a scaling that is suitable for every regularity level $\alpha$ simultaneously. Furthermore, we extend the earlier results to multivariate functions, and show that the procedure also adapts to a scale of infinitely smooth functions, of the type considered in [4, 20, 22, 23] and [32]. The proofs of several lemmas have common elements with [47], but the main result is proved from first principles.

Of course, a (rescaled) Gaussian random field is not a suitable model for a density or a binary regression function. Following other authors we transform it for these applications by exponentiation and renormalization, or by application of a link function. These transformations and the statistical consequences for these settings are given in Section 2, together with the application to the regression model. In Section 3 we state a more abstract result on rescaled Gaussian random fields, which gives the common structure to the three statistical applications. This abstract result also applies to other statistical settings, not discussed in this paper, and concerns Gaussian random fields more general than the squared exponential process, and bandwidths more general than the inverse Gamma. Proofs are deferred to Sections 4 and 5.

We consider only compactly supported functions as parameters, even though the priors in principle are functions on the full Euclidean space. Consistency of a posterior on the full space can be expected only if the tails of the functions are restricted. If they are not, then one would still expect that the posterior restricted to compact subsets contracts at some rate. At the moment there seem to exist no results that would yield such a rate (or even consistency).

1.1. *Notation.* Let $C[0,1]^d$ and $C^\alpha[0,1]^d$ be the space of all continuous functions and the Hölder space of $\alpha$-smooth functions $f:[0,1]^d \to \mathbb{R}$, respectively, equipped with the uniform norm $\|\cdot\|_\infty$ (cf. [45], Section 2.7.1). Let $\mathcal{A}^{\gamma,r}(\mathbb{R}^d)$ be the space of functions $f:\mathbb{R}^d \to \mathbb{R}$ with Fourier transform $\hat{f}(\lambda) = (2\pi)^{-d}\int e^{i(\lambda,t)}f(t)\,dt$ satisfying $\int e^{\gamma\|\lambda\|^r}|\hat{f}|^2(\lambda)\,d\lambda < \infty$. These functions are infinitely often differentiable and "increasingly smooth" as $\gamma$ or $r$ increase; they extend to functions that are analytic on a strip in $\mathbb{C}^d$ containing $\mathbb{R}^d$ if $r=1$ and to entire functions if $r>1$ (see, e.g., [3], Theorem 8.3.5).

**2. Main results.** In this section we present the main results for three different statistical settings: i.i.d. density estimation, fixed design regression and classification. The proofs of these results are consequences of a theorem on rescaled Gaussian processes in Section 3, general posterior convergence rate results from [16] and [17] and results mapping the three settings to these general results given in [46]. The process $W$ and variable $A^d$ in this section are taken to be the squared exponential Gaussian field and an independent random variable with a Gamma distribution. For $W$ and $A$ satisfying the more general conditions given in Section 3, the same results are true, except for the fact that the powers of the logarithmic factors may be different.

2.1. *Density estimation.* After exponentiation and renormalization a randomly rescaled Gaussian process can be used as a prior model for probability densities. Priors of this type were, among others, considered by [29, 30] and [31]. Posterior consistency was recently obtained in the paper [43].

To describe our adaptation result, consider a sample $X_1,\ldots,X_n$ from a continuous, positive density $f_0$ on the unit cube $[0,1]^d \subset \mathbb{R}^d$. As a prior distribution $\Pi$ on $f_0$ we use the distribution of

$$(2.1) \qquad t \mapsto \frac{e^{W_{At}}}{\int_{[0,1]^d} e^{W_{As}}\,ds}.$$

Let $\Pi(f \in \cdot | X_1, \ldots, X_n)$ denote the posterior distribution: the conditional distribution of $f$ on the Borel sets in $C[0,1]^d$ in the Bayesian setup, where the density $f$ is first drawn from the prior (2.1) and given $f$ the variables are an i.i.d. sample from $f$. We say that the posterior contracts at rate $\varepsilon_n$ if, for every sufficiently large constant $M$, as $n \to \infty$,

$$\Pi(f : h(f, f_0) \geq M\varepsilon_n | X_1, \ldots, X_n) \xrightarrow{P_{f_0}} 0.$$

Here $h$ is the Hellinger distance and the convergence is understood to be in probability under the (frequentist) assumption that $X_1, \ldots, X_n$ are a random sample from $f_0$.

THEOREM 2.1. *Let $w_0 = \log f_0$.*

- *If $w_0 \in C^\alpha[0,1]^d$ for some $\alpha > 0$, then the posterior contracts at rate $n^{-\alpha/(2\alpha+d)}(\log n)^{(4\alpha+d)/(4\alpha+2d)}$.*
- *If $w_0$ is the restriction of a function in $\mathcal{A}^{\gamma,r}(\mathbb{R}^d)$, then the posterior contracts at rate $n^{-1/2}(\log n)^{d+1}$ if $r \geq 2$ and $n^{-1/2}(\log n)^{d+1+d/(2r)}$ if $r < 2$.*

The minimax rate of estimation of a density $f_0$ that is bounded away from zero and known to belong to the space $C^\alpha[0,1]^d$ of $\alpha$-Hölder continuous functions is $n^{-\alpha/(2\alpha+d)}$. The first assertion of the theorem shows that the posterior contracts at the minimax rate times a logarithmic factor. It is rate-adaptive in the sense that this is true for any $\alpha > 0$, even though the prior does not depend on $\alpha$. We conjecture that a logaritmic factor in the rate for the present prior is necessary, although the power $(4\alpha+d)/(4\alpha+2d)$ may not be optimal. As shown in Section 3 this power can be improved by using a slightly different prior for $A$. Other Bayesian schemes (see, e.g., [18, 21] and [28]) give adaptation without logarithmic factors, but are more complicated.

The second assertion shows that the rate improves to $1/\sqrt{n}$ times a logarithmic factor if $\log f_0$ is the restriction of a function in $\mathcal{A}^{\gamma,r}(\mathbb{R}^d)$. The rate is better if $r$ increases, but does not improve beyond $r = 2$, the exponent of the spectral density of the squared exponential process. For a Gaussian prior with a compactly supported spectral density, the rate would strictly improve as $r$ increases, reaching the rate $n^{-1/2}(\log n)^{d+1}$ as $r \uparrow \infty$. Other estimation schemes (see [4, 20, 22, 23] and [32]) can reach the better rate $n^{-1/2}(\log n)^{(d+1)/2}$.

2.2. *Fixed design regression.* Suppose we observe independent variables $Y_1, \ldots, Y_n$ satisfying the regression relation $Y_i = w_0(t_i) + \varepsilon_i$, for independent $N(0, \sigma_0^2)$-distributed error variables $\varepsilon_i$ and known elements $t_1, \ldots, t_n$ of the unit cube $[0,1]^d$. The aim is to estimate the regression function $w_0$. In this case a rescaled Gaussian process can be used directly as a prior for $w_0$; cf.

[24, 49] and [40]. Posterior consistency for priors of this type was recently established in [9].

We use law of the random field $(W_{At}: t \in [0,1]^d)$ as a prior for $w_0$. If the standard deviation $\sigma_0$ of the errors is unknown, we endow it with a prior distribution as well, which we assume to be supported on a given interval $[a,b] \subset (0,\infty)$ that contains $\sigma_0$, with a Lebesgue density that is bounded away from zero.

We denote the posterior distribution by $\Pi(\cdot|Y_1,\ldots,Y_n)$. Let $\|w\|_n = (n^{-1} \times \sum_{i=1}^n w^2(t_i))^{1/2}$ be the $L_2$-norm corresponding to the empirical distribution of the design points. We say that the posterior contracts at rate $\varepsilon_n$ if, for every sufficiently large $M$,

$$\Pi((w,\sigma): \|w-w_0\|_n + |\sigma - \sigma_0| \geq M\varepsilon_n | Y_1, \ldots, Y_n) \stackrel{P_{(w_0,\sigma_0)}}{\longrightarrow} 0.$$

THEOREM 2.2. *The assertions of Theorem 2.1 are true in the setting of regression for $w_0 = f_0$.*

2.3. *Classification.* In the setting of classification, or binary regression, the use of rescaled Gaussian process priors was considered for instance in [7] and [40]. Consistency results were obtained in [14] and more recently in [8].

Consider i.i.d. observations $(X_1,Y_1),\ldots,(X_n,Y_n)$, where $X_i$ takes values in the unit cube $[0,1]^d$ and $Y_i$ takes values in the set $\{0,1\}$. The statistical problem is to estimate the binary regression function $r_0(t) = P(Y_1 = 1|X_1 = t)$.

As a prior $\Pi$ on $r_0$ we use the law of the process $(\Psi(W_{At}): t \in [0,1]^d)$, where $\Psi: \mathbb{R} \to (0,1)$ is the logistic or the normal distribution function.

Let $\Pi(\cdot|(X_1,Y_1),\ldots,(X_n,Y_n))$ denote the posterior and let $\|\cdot\|_{L_2(G)}$ be the $L_2$-norm relative to the marginal distribution $G$ of $X_1$. We say that the posterior contracts at rate $\varepsilon_n$ if, for every sufficiently large $M$,

$$\Pi(r: \|r - r_0\|_{L_2(G)} \geq M\varepsilon_n|(X_1,Y_1),\ldots,(X_n,Y_n)) \stackrel{P_{r_0}}{\longrightarrow} 0.$$

THEOREM 2.3. *Let $w_0 = \Psi^{-1}(r_0)$. Then the assertions of Theorem 2.1 are true.*

**3. Rescaled Gaussian fields.** Let $W = (W_t: t \in \mathbb{R}^d)$ be a centered, homogeneous Gaussian random field with covariance function of the form, for a given continuous function $\phi: \mathbb{R}^d \to \mathbb{R}$,

$$(3.1) \qquad EW_sW_t = \phi(s-t).$$

By Bochner's theorem there exists a finite Borel measure $\mu$ on $\mathbb{R}^d$, the *spectral measure* of $W$, such that

$$\phi(t) = \int e^{-i(\lambda,t)} \mu(d\lambda). \tag{3.2}$$

We shall consider processes whose spectral measure $\mu$ has subexponential tails: for some $\delta > 0$,

$$\int e^{\delta \|\lambda\|} \mu(d\lambda) < \infty. \tag{3.3}$$

The squared exponential process, whose covariance function is given in (1.1), falls in this class. Its spectral measure has density relative to the Lebesgue measure given by $\lambda \mapsto \exp(-\|\lambda\|^2/4)/(2^d \pi^{d/2})$.

For a positive random variable $A$ defined on the same probability space as $W$ and stochastically independent of $W$ let $W^A = (W_{At} : t \in [0,1]^d)$ be the restriction to $[0,1]^d$ of the rescaled process $t \mapsto W_{At}$. We consider it as a Borel measurable map in the space $C[0,1]^d$, equipped with the uniform norm $\|\cdot\|_\infty$. The following theorem bounds the small-ball probability and the complexity of the support of the field $W^A$. These are the essential ingredients for proving the statistical results in Section 2, and can also be used to analyse other Bayesian schemes.

We assume that the distribution of $A$ possesses a Lebesgue density $g$ satisfying, for positive constants $C_1, D_1, C_2, D_2$, nonnegative constants $p, q$, and every sufficiently large $a > 0$,

$$C_1 a^p \exp(-D_1 a^d \log^q a) \leq g(a) \leq C_2 a^p \exp(-D_2 a^d \log^q a). \tag{3.4}$$

This is satisfied (with $q = 0$) if $A^d$ possesses a Gamma distribution.

For given sequences $\varepsilon_n$ and $\bar{\varepsilon}_n$ and a given function $w_0 : [0,1]^d \to \mathbb{R}$, consider the following statement: there exist Borel measurable subsets $B_n$ of $C[0,1]^d$ and a constant $K$ such that, for every sufficiently large $n$,

$$P(\|W^A - w_0\|_\infty \leq \varepsilon_n) \geq e^{-n\varepsilon_n^2}, \tag{3.5}$$

$$P(W^A \notin B_n) \leq e^{-4n\varepsilon_n^2}, \tag{3.6}$$

$$\log N(\bar{\varepsilon}_n, B_n, \|\cdot\|_\infty) \leq n\bar{\varepsilon}_n^2. \tag{3.7}$$

THEOREM 3.1. *Let $W$ be a centered homogeneous Gaussian field with spectral measure $\mu$ that satisfies (3.3) for some $\delta > 0$ and that possesses a Lebesgue density $f$ such that $a \mapsto f(a\lambda)$ is decreasing on $(0,\infty)$ for every $\lambda \in \mathbb{R}^d$:*

- *If $w_0 \in C^\alpha[0,1]^d$ for some $\alpha > 0$, then there exist Borel measurable subsets $B_n$ of $C[0,1]^d$ such that (3.5), (3.6) and (3.7) hold, for every sufficiently large $n$, and $\varepsilon_n = n^{-\alpha/(2\alpha+d)}(\log n)^{\kappa_1}$ for and $\bar{\varepsilon}_n = K\varepsilon_n(\log n)^{\kappa_2}$, for $\kappa_1 = ((1+d) \vee q)/(2+d/\alpha)$ and $\kappa_2 = (1+d-q)/2$ and a sufficiently large constant $K$.*
- *If $w_0$ is the restriction of a function in $\mathcal{A}^{\gamma,r}(\mathbb{R}^d)$ to $[0,1]^d$ and the spectral density satisfies $|f(\lambda)| \geq C_3 \exp(-D_3\|\lambda\|^\nu)$ for some positive constants $C_3$, $D_3$ and $\nu$, then there exist Borel measurable subsets $B_n$ of $C[0,1]^d$ such that (3.5), (3.6) and (3.7) hold, for every sufficiently large $n$, and $\varepsilon_n = Kn^{-1/2}(\log n)^{(d+1)/2}$ for $r \geq \nu$, $\varepsilon_n = Kn^{-1/2}(\log n)^{(d+1)/2+d/(2r)}$ for $r < \nu$, and $\bar{\varepsilon}_n = \varepsilon_n(\log n)^{(d+1)/2}$, for a sufficiently large constant $K$.*

In the paper [46] it is shown that (3.5)–(3.7) map one-to-one to the general conditions on rates of contraction of posterior distributions used in [17] and [16], for each of the three settings considered in Section 2. Thus a rate of contraction $\varepsilon_n \vee \bar{\varepsilon}_n$ is attained for each of these three settings. Theorems 2.1–2.3 follow, with the parameter $q$ equal to 0. (The use of two rates $\varepsilon_n$ and $\bar{\varepsilon}_n$ requires a slight generalization of the main result in [17], formulated as Theorem 2.1 in [15]; also see the discussion following the statement of the main result in [16].) The choice $q = d+1$ yields a slightly better rate (a lower power on the logarithmic factor), but we highlighted the choice $q = 0$ in Section 2, as this corresponds to a Gamma prior.

**4. Auxiliary results.** In this section we prepare a number of auxiliary lemmas needed in the proof of Theorem 3.1. In the proof of (3.5) we condition on the variable $A$, so that we can first consider the probability in (3.5) for $A$ a fixed constant, and then combine the obtained bound with bounds on the tails of the distribution of $A$. The proofs of (3.6) and (3.7) involve similar steps.

For fixed $A$ the process $W^A$ is a Gaussian random field with values in $C[0,1]^d$, and a key concept is the associated *reproducing kernel Hilbert space* (RKHS). This can be viewed as a subset of the space $C[0,1]^d$, which gives the "geometry" of the distribution of $W^A$, just as finite-dimensional Gaussian vectors are described by ellipsoids. According to general Gaussian process theory, obtaining good bounds for the probabilities in (3.5) and (3.6) for fixed $A$ is closely linked to studying the metric entropy of the unit ball of the RKHS and the approximation of the function $w_0$ by elements of the RKHS. See [48] for a review.

In Lemma 4.1 we start by characterizing the RKHS of the process $W$, from which the RKHS of the rescaled process $W^A$ will be obtained in Lemma 4.2. The RKHS of a Gaussian field $(W_t : t \in T)$, with parameter set equal to a set $T \subset \mathbb{R}^d$, is by definition the set of functions $h : T \to \mathbb{R}$ that can be represented as $h(t) = EW_t L$ for $L$ contained in the closure of the linear span of the variables $(W_t : t \in T)$ in $L_2(\Omega, \mathcal{U}, P)$, for $(\Omega, \mathcal{U}, P)$ the probability space on which $W$ is defined, equipped with the square norm $\|h\|_{\mathbb{H}}^2 = EL^2$.

LEMMA 4.1. *The RKHS of $(W_t : t \in T)$ is the set of real parts of the functions (from $T$ to $\mathbb{C}$)*

$$t \mapsto \int e^{i(\lambda,t)} \psi(\lambda) \mu(d\lambda),$$

*when $\psi$ runs through the complex Hilbert space $L_2(\mu)$. The RKHS-norm of the displayed function equals the norm in $L_2(\mu)$ of the projection of $\psi$ on the closed linear span of the set of functions $(e_s : s \in T)$ (or, equivalently, the infimum of $\|\psi\|_2$ over all functions $\psi$ giving the same function in the preceding display). If $T \subset \mathbb{R}^d$ has an interior point and (3.3) holds, then this closed linear span is $L_2(\mu)$ and the RKHS norm is $\|\psi\|_{L_2(\mu)}$.*

PROOF. The spectral representation (3.2) can be written as $EW_s W_t = \langle e_t, e_s \rangle_{L_2(\mu)}$, for $e_t$ the function defined by $e_t(\lambda) = \exp(i(\lambda, t))$. By definition the RKHS is therefore the set of functions as in the display, with $\psi$ running through the closure $\mathbb{L}_T$ in $L_2(\mu)$ of the linear span of the set of functions $(e_s : s \in T)$, and the norm equal to the norm of $\psi$ in $L_2(\mu)$. Here the "linear span" is taken over the reals. If instead we take the linear span over the complex numbers, we obtain complex functions whose real parts give the RKHS.

The set of functions obtained by letting $\psi$ range the full space $L_2(\mu)$ is precisely the same, as a general element $\psi \in L_2(\mu)$ gives exactly the same function as its projection $\Pi \psi$ on $\mathbb{L}_T$. However, the associated norm is the $L_2(\mu)$ norm of $\Pi \psi$. This proves the first assertion of the lemma. For the second we must show that $\mathbb{L}_T = L_2(\mu)$ under the additional conditions.

The partial derivative of order $(k_1, \ldots, k_d)$ with respect to $(t_1, \ldots, t_d)$ of the map $t \mapsto e_t$ at $t_0$ is the function $\lambda \mapsto (i\lambda_1)^{k_1} \cdots (i\lambda_d)^{k_d} e_{t_0}(\lambda)$. Appealing to the dominated convergence theorem we see that this derivative exists as a derivative in $L_2(\mu)$. Because $t_0$ is an interior point of $T$ by assumption, we conclude that the function $\lambda \mapsto (i\lambda)^k e_{t_0}(\lambda)$ belongs to $\mathbb{L}_T$ for any multiindex $k$ of nonnegative integers. Consequently, the function $p e_{t_0}$ belongs to $\mathbb{L}_T$ for any polynomial $p : \mathbb{R}^d \to \mathbb{C}$ in $d$ arguments. It suffices to show that these functions are dense in $L_2(\mu)$.

Equivalently, it suffices to prove that the polynomials themselves are dense in $L_2(\mu)$. Indeed, if $\psi \in L_2(\mu)$ is orthogonal to all functions of the form $p e_{t_0}$, then $\psi \overline{e_{t_0}}$ is orthogonal to all polynomials. Denseness of the set of polynomials then gives that $\psi \overline{e_{t_0}}$ vanishes $\mu$-almost everywhere, whence $\psi$ vanishes $\mu$-almost everywhere.

That the polynomials are dense in $L_2(\mu)$ appears to be well known. A proof for $d = 1$ is given in [38]. For completeness we include a proof for general dimension $d$. Suppose that $\psi \in L_2(\mu)$ is orthogonal to all polynomials.

Since $\mu$ is a finite measure, the complex conjugate $\overline{\psi}$ is $\mu$-integrable, and hence we can define a complex measure $\nu$ by

$$\nu(B) = \int_B \overline{\psi(\lambda)} \mu(d\lambda).$$

It suffices to show that $\nu$ is the zero measure, so that $\psi = 0$ almost everywhere relative to $\mu$.

By the Cauchy–Schwarz inequality and (3.3), with $\|\nu\|$ the (total) variation measure of $\nu$,

(4.1) $$\int e^{\delta\|\lambda\|/2} \|\nu\|(d\lambda) < \infty.$$

By a standard argument, based on the dominated convergence theorem (see, e.g., [3], Theorem 8.3.5), this implies that the function $z \mapsto \int e^{(\lambda,z)} \nu(d\lambda)$ is analytic on the strip $\Omega = \{z \in \mathbb{C}^d : |\operatorname{Re} z_1| < \delta/(2\sqrt{d}), \ldots, |\operatorname{Re} z_d| < \delta/(2\sqrt{d})\}$. Also for $z$ real and in this strip, by the dominated convergence theorem,

$$\int e^{(\lambda,z)} \nu(d\lambda) = \int \sum_{n=0}^{\infty} \frac{(\lambda,z)^n}{n!} \nu(d\lambda) = \sum_{n=0}^{\infty} \int \frac{(\lambda,z)^n}{n!} \overline{\psi}(\lambda) \mu(d\lambda).$$

The right-hand side vanishes, because $\psi$ is orthogonal to all polynomials by assumption.

We conclude that the function $z \mapsto \int e^{(\lambda,z)} \nu(d\lambda)$ vanishes on the set $\{z \in \Omega : \operatorname{Im} z = 0\}$. Because this set contains a nontrivial interval in $\mathbb{R}$ for every coordinate, we can apply (repeated) analytic continuation to see that this function vanishes on the complete strip $\Omega$. In particular the Fourier transform $t \mapsto \int e^{i(\lambda,t)} \nu(d\lambda)$ of $\nu$ vanishes on all of $\mathbb{R}^d$, whence $\nu$ is the zero-measure. □

For $W = (W_t : t \in \mathbb{R}^d)$ a homogeneous Gaussian random field with spectral measure $\mu$ and a positive real number $a$, the rescaled process $(W_{at} : t \in \mathbb{R}^d)$ is also homogenous and has spectral measure $\mu_a$ that is related to $\mu$ by

$$\mu_a(B) = \mu(B/a).$$

If $\mu$ has a (spectral) density $f$, then $\mu_a$ has density $f_a$ given by

$$f_a(\lambda) = a^{-d} f(\lambda/a).$$

We shall obtain approximation properties and small-ball probabilities for the process $W^a = (W_{at} : t \in [0,1]^d)$, viewed as a map in $C[0,1]^d$. Let $\mathbb{H}^a$ be the RKHS of $W^a$, with corresponding norm $\|\cdot\|_{\mathbb{H}^a}$. It is described in Lemma 4.1 with $\mu$ taken equal to $\mu_a$.

The following lemma follows from general principles, or can be proved from the characterization of RKHSs given in Lemma 4.1. By "scaling map" $h \mapsto (t \mapsto h(at))$ we mean the map that attaches to a given function $h : [0,a]^d \to \mathbb{R}$ the function $g : [0,1]^d \to \mathbb{R}$ defined by $g(t) = h(at)$.

LEMMA 4.2. *The scaling map $h \mapsto (t \mapsto h(at))$ is an isometry from the RKHS of the process $(W_t : t \in [0,a]^d)$ onto $\mathbb{H}^a$.*

The next step is to bound the concentration function of the Gaussian prior $W^a$, again for a fixed $a$. The concentration function (at $\varepsilon > 0$) is the sum of minus the log centered small probability, considered in Lemma 4.6, and the decentering function $\inf\{\|h\|_{\mathbb{H}^a}^2 : \|h - w_0\|_\infty < \varepsilon\}$, which measures the positioning of the true parameter $w_0$ relative to the RKHS. We start by bounding the latter, separately for the cases that the true parameter is Hölder or supersmooth in Lemmas 4.3 and 4.4. The first lemma is fairly standard, and proceeds by approximating $w_0$ by a suitable convolution of $w_0$ with a smooth function, which is contained in the RKHS.

LEMMA 4.3. *Assume that the restriction of $\mu$ to some neighborhood of the origin is Lebesgue absolutely continuous with a density that is bounded away from zero. Let $\alpha > 0$ be given. Then for any $w \in C^\alpha[0,1]^d$ there exist constants $C$ and $D$ depending only on $\mu$ and $w$ such that, as $a \to \infty$,*

$$\inf\{\|h\|_{\mathbb{H}^a}^2 : \|h - w\|_\infty \leq Ca^{-\alpha}\} \leq Da^d.$$

PROOF. Let $\underline{\alpha}$ be the biggest integer strictly smaller than $\alpha$. Let $G$ be a bounded neighborhood of the origin on which $\mu$ has a Lebesgue density $f$ that is bounded away from 0. Take a function $\psi : \mathbb{R} \to \mathbb{C}$ with a symmetric, real-valued, infinitely smooth Fourier transform $\hat\psi$ that is supported on an interval $I$ such that $I^d \subset G$ and which equals $1/(2\pi)$ in a neighborhood of zero, so that $\psi$ has moments of all orders and

$$\int (it)^k \psi(t)\, dt = 2\pi \hat\psi^{(k)}(0) = \begin{cases} 0, & k \geq 1, \\ 1, & k = 0. \end{cases}$$

Define $\phi : \mathbb{R}^d \to \mathbb{C}$ by $\phi(t) = \psi(t_1)\cdots\psi(t_d)$. Then we have that $\int \phi(t)\, dt = 1$, and $\int t^k \phi(t)\, dt = 0$, for any nonzero multi-index $k = (k_1, \ldots, k_d)$ of nonnegative integers. Moreover, we have that $\int \|t\|^\alpha |\phi|(t)\, dt < \infty$, and the functions $|\hat\phi|/f$ and $|\hat\phi|^2/f$ are uniformly bounded.

By Whitney's theorem we can extend $w : [0,1]^d \to \mathbb{R}$ to a function $w : \mathbb{R}^d \to \mathbb{R}$ with compact support and $\|w\|_\alpha < \infty$. (See [50] or [41], Chapter VI; we can multiply an arbitrary smooth extension by an infinitely smooth function that vanishes outside a neighborhood of $[0,1]^d$ to ensure compact support).

By Taylor's theorem we can write, for $s, t \in \mathbb{R}^d$,

$$w(t+s) = \sum_{j\,:\,j.\leq\underline\alpha} D^j w(t) \frac{s^j}{j!} + S(t,s),$$

where

$$|S(t,s)| \leq C\|s\|^\alpha$$

for a positive constant $C$ that depends on $w$ but not on $s$ and $t$. If we set $\phi_a(t) = \phi(at)$ we get, in view of the fact that $\phi$ is a higher-order kernel, for any $t \in \mathbb{R}^d$,

$$a^d(\phi_a * w)(t) - w(t) = \int \phi(s)(w(t - s/a) - w(t))\,ds = \int \phi(s) S(t, -s/a)\,ds.$$

Combining the preceding displays shows that $\|a^d \phi_a * w - w\|_\infty \leq KCa^{-\alpha}$, for $K = \int \|s\|^\alpha |\phi|(s)\,ds$.

For $\hat{w}$ the Fourier transform of $w$, we can write

$$\frac{1}{(2\pi)^d}(\phi_a * w)(t) = \int e^{-i(t,\lambda)} \hat{w}(\lambda) \hat{\phi}_a(\lambda)\,d\lambda = \int e^{-i(t,\lambda)} \frac{\hat{w}(-\lambda)\hat{\phi}_a(\lambda)}{f_a(\lambda)}\,d\mu_a(\lambda).$$

Therefore, by Lemma 4.1 the function $a^d \phi_a * w$ is contained in the RKHS $\mathbb{H}^a$, with square norm a multiple of, with $\Pi$ the orthogonal projection in $L_2(\mu)$ onto the functions $(e_t : t \in [0,1]^d)$,

$$a^{2d} \int \left|\Pi\left(\frac{\hat{w}\hat{\phi}_a}{f_a}\right)\right|^2 d\mu_a \leq a^d \int \frac{|\hat{w}(\lambda)|^2 |\hat{\phi}(\lambda/a)|^2}{f(\lambda/a)}\,d\lambda$$

$$\leq a^d \int |\hat{w}(\lambda)|^2\,d\lambda \left\|\frac{|\hat{\phi}|^2}{f}\right\|_\infty.$$

Here $(2\pi)^d \int |\hat{w}|^2(\lambda)\,d\lambda = \int |w|^2(t)\,dt$ is finite, and $|\hat{\phi}|^2/f$ is bounded by the construction of $\hat{\phi}$. $\square$

The supersmooth case consists of the subcase that $w_0$ is "super-super smooth," that is, it belongs itself to the RKHS, and the more regular case in which it is approximated by its "projection" in the RKHS.

LEMMA 4.4. *Assume that $\mu$ has a Lebesgue density $f$ such that $|f(\lambda)| \geq C_3 \exp(-D_3 \|\lambda\|^\nu)$ for some positive constants $C_3$, $D_3$ and $\nu$.*

- *If $w$ is the restriction to $[0,1]^d$ of an element of $\mathcal{A}^{\gamma,r}(\mathbb{R}^d)$ for $r \geq \nu$, then $w \in \mathbb{H}^a$ for all sufficiently large $a$ with uniformly bounded norm $\|w\|_{\mathbb{H}^a}$.*
- *If $w$ is the restriction to $[0,1]^d$ of an element of $\mathcal{A}^{\gamma,r}(\mathbb{R}^d)$ for $r < \nu$, then there exist constants $a_0$, $C$ and $D$ depending only on $\mu$ and $w$ such that, for $a > a_0$,*

$$\inf\{\|h\|_{\mathbb{H}^a}^2 : \|h - w\|_\infty \leq Ce^{-\gamma a^r}/a^{-r+1}\} \leq Da^d.$$

PROOF. The Fourier transform of a function $w \in \mathcal{A}^{\gamma,r}(\mathbb{R}^d)$ is certainly integrable, and hence, by the inversion formula,

$$w(t) = \int e^{-i(\lambda,t)} \hat{w}(\lambda)\,d\lambda = \int e^{-i(\lambda,t)} \frac{\hat{w}}{f_a}(\lambda)\,d\mu_a(\lambda).$$

In view of Lemma 4.1 $w \in \mathbb{H}^a$ if $\hat{w}/f_a \in L_2(\mu_a)$. Now

$$\int \left|\frac{\hat{w}}{f_a}\right|^2 d\mu_a \leq \int |\hat{w}(\lambda)|^2 \frac{a^d}{C_3} e^{D_3 \|\lambda\|^\nu / a^\nu} \, d\lambda.$$

This is finite for every $a > 0$ if $r > \nu$. If $r = \nu$, then this is finite for $a \geq (D_3/\gamma)^{1/\nu}$. In both cases the right side is bounded as $a \to \infty$.

To prove the second assertion let $\phi$ be as in the proof of Lemma 4.3, with compactly supported Fourier transform $\hat{\phi}$ constructed to be constant and equal to $(2\pi)^{-d}$ on $[-1,1]^d$, and bounded in absolute value by this constant everywhere. By the argument given in this proof the function $a^d \phi_a * w$ is contained in $\mathbb{H}^a$ with square norm bounded above by a multiple of $a^d$, for sufficiently large $a$. Also

$$|a^d \phi_a * w(t) - w(t)|^2 = \left|\int e^{-i(\lambda,t)} \left((2\pi)^d \hat{\phi}\left(\frac{\lambda}{a}\right) - 1\right) \hat{w}(\lambda) \, d\lambda\right|^2$$
$$\leq \left(\int_{\|\lambda/a\| > 1} 2|\hat{w}(\lambda)| \, d\lambda\right)^2$$
$$\leq 4 \int_{\|\lambda\| > a} e^{-\gamma \|\lambda\|^r} \, d\lambda \int |\hat{w}(\lambda)|^2 e^{\gamma \|\lambda\|^r} \, d\lambda,$$

by the Cauchy–Schwarz inequality. The second factor is finite if $w \in \mathcal{A}^{\gamma,r}(\mathbb{R}^d)$. The first is bounded by a multiple of $e^{-\gamma a^r} a^{-r+1}$, by a change of variable and Lemma 4.9. □

Next we turn to bounding the centered small-ball probability. According to general results on Gaussian processes (see [26]), this can be characterized in terms of the entropy of the unit ball of the RKHS. In view of Lemma 4.1 this consists of certain analytic functions, and therefore we can bound its entropy by employing classical techniques as given in [25].

Let $\mathbb{H}_1^a$ be the unit ball in the RKHS of $W^a = (W^a : t \in [0,1]^d)$, that is, the set of functions $h \in \mathbb{H}^a$ with $\|h\|_{\mathbb{H}^a} \leq 1$.

LEMMA 4.5.  *Let $\mu$ satisfy (3.3) for some $\delta > 0$. There exists a constant $K$, depending only on $\mu$ and $d$, such that, for $\varepsilon < 1/2$,*

$$\log N(\varepsilon, \mathbb{H}_1^a, \|\cdot\|_\infty) \leq K a^d \left(\log \frac{1}{\varepsilon}\right)^{1+d}.$$

PROOF.  By Lemma 4.1 a typical element of $\mathbb{H}_1^a$ can be written as the real part of the function $h_\psi : [0,1]^d \to \mathbb{C}$ given by

(4.2)  $$h_\psi(t) = \int e^{i(\lambda,t)} \psi(\lambda) \mu_a(d\lambda),$$

for $\psi\colon\mathbb{R}^d\to\mathbb{C}$ a function with $\int|\psi|^2\mu_a(d\lambda)\le 1$. We shall construct an $\varepsilon$-net over these functions consisting of piecewise polynomials.

For $R=\delta/(3a\sqrt{d})$ let $\{t_1,\ldots,t_m\}$ be an $R/2$-net in $T=[0,1]^d$, for the maximum norm, and let $T=\bigcup_i B_i$ be a partition of $T$ in sets $B_1,\ldots,B_m$ obtaining by assigning every $t\in T$ to a closest $t_i\in\{t_1,\ldots,t_m\}$. Consider the piecewise polynomials $P=\sum_{i=1}^m P_{i,a_i}1_{B_i}$, for

$$P_{i,a_i}(t)=\sum_{n.\le k} a_{i,n}(t-t_i)^n.$$

Here the sum ranges over all multi-index vectors $n=(n_1,\ldots,n_d)\in(\mathbb{N}\cup\{0\})^d$ with $n.=n_1+\cdots+n_d\le k$, and for $s=(s_1,\ldots,s_d)\in\mathbb{R}^d$ the notation $s^n$ is short for $s_1^{n_1}s_2^{n_2}\cdots s_d^{n_d}$. We obtain a finite set of functions by discretizing the coefficients $a_{i,n}$ for each $i$ and $n$ over a grid of meshwidth $\varepsilon/R^{n.}$-net in the interval $[-C/R^{n.},C/R^{n.}]$, for given $C>0$. The log cardinality of this set is bounded by

$$\log\Bigl(\prod_i\prod_{n:n.\le k}\#a_{i,n}\Bigr)\le m\log\Bigl(\prod_{n:n.\le k}\frac{2C/R^{n.}}{\varepsilon/R^{n.}}\Bigr)\le mk^d\log\Bigl(\frac{2C}{\varepsilon}\Bigr).$$

We can choose $m\le(3/R/2)^d$. The proof is complete once it is shown that the resulting set of functions is a $K\varepsilon$-net for constants $C$ and $K$ depending only on $\mu$, and for $k$ of the order $\log(1/\varepsilon)$.

We can view the function $h_\psi$ as a function of the argument $it$, ranging over the product of the imaginary axes in $\mathbb{C}^d$. In view of (3.3) and the Cauchy–Schwarz inequality, this function can be extended to an analytic function $z\mapsto\int e^{(\lambda,z)}\psi(\lambda)\,d\mu_a(\lambda)$ on the set $\{z\in\mathbb{C}^d\colon\|\operatorname{Re}z\|<\delta/2\}$, which includes the strip $\Omega=\{z\in\mathbb{C}^d\colon|\operatorname{Re}z_1|\le R,\ldots,|\operatorname{Re}z_d|\le R\}$ for $R=\delta/(3a\sqrt{d})$, and it satisfies the uniform bound, for every $z\in\Omega$,

$$|h_\psi(z)|^2\le\int e^{\delta\|\lambda\|}\mu(d\lambda):=C^2.$$

By the Cauchy formula ($d$ applications of the formula in one dimension suffice), for $C_1,\ldots,C_d$ circles of radius $R$ in the complex plane around the coordinates $t_{i1},\ldots,t_{id}$ of $t_i$, and with $D^n$ the partial derivative of orders $n=(n_1,\ldots,n_d)$ and $n!=n_1!n_2!\cdots n_d!$,

$$\Bigl|\frac{D^n h_\psi(t_i)}{n!}\Bigr|=\Bigl|\frac{1}{(2\pi i)^d}\oint_{C_1}\cdots\oint_{C_d}\frac{h_\psi(z)}{(z-t_i)^{n+1}}\,dz_1\cdots dz_d\Bigr|\le\frac{C}{R^{n.}}.$$

Consequently, for any $z\in B_i$, a universal constant $K$, and appropriately chosen $a_i$

$$\Bigl|\sum_{n.>k}\frac{D^n h_\psi(t_i)}{n!}(z-t_i)^n\Bigr|\le\sum_{n.>k}\frac{C}{R^{n.}}(R/2)^{n.}\le C\sum_{l=k+1}^\infty\frac{l^{d-1}}{2^l}$$

$$\leq KC\left(\frac{2}{3}\right)^k,$$

$$\left|\sum_{n.\leq k}\frac{D^n h_\psi(t_i)}{n!}(z-t_i)^n - P_{i,a_i}(z)\right| \leq \sum_{n.\leq k}\frac{\varepsilon}{R^{n.}}(R/2)^{n.} \leq \sum_{l=1}^{k}\frac{l^{d-1}}{2^l}\varepsilon \leq K\varepsilon.$$

We conclude that the piecewise polynomials form a $2K\varepsilon$-net for $k$ sufficiently large that $(2/3)^k$ is smaller than $K\varepsilon$. □

LEMMA 4.6. *If the spectral measure satisfies (3.3), then for any $a_0 > 0$ there exists constants $C$ and $\varepsilon_0$ that depend only on $a_0$, $\mu$ and $d$ only such that, for $a \geq a_0$ and $\varepsilon < \varepsilon_0$,*

$$-\log P\Big(\sup_{t\in[0,1]^d}|W_t^a| \leq \varepsilon\Big) \leq Ca^d\left(\log\frac{a}{\varepsilon}\right)^{1+d}.$$

PROOF. This is essentially a corollary of Lemma 4.5 in the present paper and Theorem 2 of [26]. However, to make the dependence on the scaling factor $a$ explicit it is necessary to go through the steps of the proof of the latter theorem. We only sketch the main steps of the long derivation. Let $\phi_0^a(\varepsilon)$ be the left side of the lemma.

By formula (3.19) of [26], for any $\varepsilon, \lambda > 0$,

$$\phi_0^a(2\varepsilon) + \log\Phi(\lambda + \Phi^{-1}(e^{-\phi_0^a(\varepsilon)})) \leq \log N\left(\frac{\varepsilon}{\lambda}, \mathbb{H}_1^a, \|\cdot\|_\infty\right).$$

Choosing $\lambda = \sqrt{2\phi_0^a(\varepsilon)}$, using the fact that $\Phi(\sqrt{2x} + \Phi^{-1}(e^{-x})) \geq 1/2$ for every $x > 0$ (see Lemma 4.10), and applying Lemma 4.5 to the right of the preceding display, we conclude that, for every $\varepsilon < 1/2$,

$$\phi_0^a(2\varepsilon) + \log\frac{1}{2} \leq Ka^d\left(\log\frac{\phi_0^a(\varepsilon)}{\varepsilon}\right)^{1+d}.$$

The (apparently) most difficult part of the proof is to show a crude bound of the form, for $\varepsilon < \varepsilon_0$ and $a \geq a_0$, and some $\tau > 0$,

$$(4.3) \qquad \phi_0^a(\varepsilon) \leq C_\tau\left(\frac{a}{\varepsilon}\right)^\tau.$$

Inserting this bound in the right of the preceding display gives that this is bounded by

$$Ka^d\left((\tau+1)\log\frac{1}{\varepsilon} + \log C_\tau + \tau\log a\right)^{1+d}.$$

This implies the assertion of the lemma.

The bound (4.3) follows for fixed $a$ immediately from Proposition 2.4 of [36], whose condition is satisfied for any $\alpha > 0$ in our case, so that we can use any $\tau > 0$. To see the dependence on $a$ we can follow the proof of Proposition 2.4, which unfortunately is involved. We only note that the constants in Lemma 2.1 of [36] (which is quoted from [39]) are universal and hence cause no problems; that Lemma 2.2 of [36] (which is quoted from [44]) can be formulated to say that $\sup_{k \leq n} k^\alpha e_k(u^*) \leq 32 \sup_{k \leq n} k^\alpha e_k(u)$ for every $n$, without conditions, and hence only involves the constant 32; finally, the proof of Proposition 2.2 is given in [36] and does not cause problems. □

For different values of $a$ the processes $W^a$ result from rescaling a single Gaussian field by different amounts. This leads to a nesting property of the attached RKHSs.

LEMMA 4.7. *Assume (3.3). If $a \leq b$, then $\sqrt{a}\mathbb{H}_1^a \subset \sqrt{b}\mathbb{H}_1^b$.*

PROOF. This follows from the characterization of the RKHS given in Lemma 4.1, together with the observations

$$\int e^{i(\lambda,t)}\psi(\lambda)\,d\mu_a(\lambda) = \int e^{i(\lambda,t)}\left(\psi\frac{f_a}{f_b}\right)(\lambda)\,d\mu_b(\lambda),$$

$$\int \left|\psi\frac{f_a}{f_b}\right|^2 d\mu_b \leq \left\|\frac{f_a}{f_b}\right\|_\infty \int |\psi|^2\,d\mu_a \leq \frac{b}{a}\int |\psi|^2\,d\mu_a.$$

Here we use that $f_a/f_b(\lambda) = (b/a)f(\lambda/a)/f(\lambda/b) \leq b/a$ by the assumed radial monotonicity of the density $f$ of the spectral measure $\mu$. □

If $a \downarrow 0$ the sample paths of $W^a$ tend on compacta to the constant value $W_0$. The following lemma gives a corresponding property for the RKHSs.

LEMMA 4.8. *Any $h \in \mathbb{H}_1^a$ satisfies $|h(0)| \leq \sqrt{\|\mu\|}$ and $|h(t) - h(0)| \leq a\|t\|\tau$ for $\tau^2 = \int \|\lambda\|^2\,d\mu(\lambda)$, for every $t \in T$.*

PROOF. By Lemma 4.1 a typical element of $\mathbb{H}_1^a$ can be written as the real part of $h(t) = \int e^{i(\lambda,t)}\psi(\lambda)\,d\mu_a(\lambda)$ for a function $\psi$ with $\int |\psi|^2\,d\mu_a \leq 1$. It follows that $|h(0)| \leq \int |\psi|\,d\mu_a$ and $|h(t) - h(0)| \leq \int |(\lambda,t)||\psi|(\lambda)\,d\mu_a(\lambda)$. Two applications of the Cauchy–Schwarz inequality conclude the proof. □

The final two lemmas in this section bound the tail probabilities of the scaling variable $A$, and give a bound on the normal quantile function, for easy reference.

LEMMA 4.9. *If the random variable $A$ has a density $g$ that satisfies (3.4) for some $q \geq 0$, then for $a^d(\log a)^q > 2|p - d + 1|/(D_2 d)$ and $a > e$,*

$$P(A > a) \leq \frac{2C_2 a^{p-d+1} \exp(-D_2 a^d (\log a)^q)}{D_2 d (\log a)^q}.$$

PROOF. Set $j_{p,r}(s) = s^p \exp(-D_2 s^d (\log s)^q)(\log s)^r$ and $J_{p,r}(a) = \int_a^\infty j_{p,r}(s)\,ds$. The derivative of the function $j_{p,0}$ can, with the help of the chain rule, be expressed as the sum of three terms. By integrating this identity we see that

$$j_{p,0}(a) = D_2 d J_{p+d-1,q}(a) + Dq J_{p,q-1}(a) - p J_{p-1,0}(a).$$

The middle term on the right is nonnegative (the third is negative if and only if $p > 0$). By the transformation $p + d - 1 \to p$ we conclude that

$$D_2 d J_{p,q}(a) - |p - d + 1| J_{p-d,0}(a) \leq j_{p-d+1,0}(a).$$

Here $J_{p,q}(a) \geq (\log a)^q J_{p,0}(a)$ and $J_{p-d,0}(a) \leq a^{-d} J_{p,0}(a)$. By substituting these inequalities in the left-hand side and rearranging we obtain the bound on $P(A > a) \leq C_2 J_{p,0}(a)$ asserted by the lemma. □

LEMMA 4.10. *The standard normal distribution function $\Phi$ satisfies $\Phi(x) \leq \exp(-x^2/2)$ for $x < 0$ and $-\sqrt{2\log(1/u)} \leq \Phi^{-1}(u)$ for $u \in (0,1)$ and $\Phi^{-1}(u) \leq -\frac{1}{2}\sqrt{\log(1/u)}$ for $u \in (0, 1/4)$.*

**5. Proof of Theorem 3.1.** For a given $a > 0$ define centered and decentered concentration functions of the process $W^a = (W_{at} : t \in [0,1]^d)$ by

$$\phi_0^a(\varepsilon) = -\log P(\|W^a\|_\infty \leq \varepsilon),$$
$$\phi_{w_0}^a(\varepsilon) = \inf_{h \in \mathbb{H}^a \,:\, \|h - w_0\|_\infty \leq \varepsilon} \|h\|_{\mathbb{H}^a}^2 - \log P(\|W^a\|_\infty \leq \varepsilon).$$

Then $P(\|W^a\|_\infty \leq \varepsilon) = \exp(-\phi_0^a(\varepsilon))$ by definition, and by results of [27] (cf. Lemma 5.3 of [48]),

(5.1) $$P(\|W^a - w_0\|_\infty \leq 2\varepsilon) \geq e^{-\phi_{w_0}^a(\varepsilon)}.$$

By Lemma 4.6 we have that $\phi_0^a(\varepsilon) \leq C_4 a^d (\log(a/\varepsilon))^{1+d}$ for $a > a_0$ and $\varepsilon < \varepsilon_0$, where the constants $a_0, \varepsilon_0, C_4$ depend only on $\mu$ and $w$.

For $\mathbb{B}_1$ the unit ball of $C[0,1]^d$ and given positive constants $M, r, \delta, \varepsilon$ set

$$B = \left(M\sqrt{\frac{r}{\delta}} \mathbb{H}_1^r + \varepsilon \mathbb{B}_1\right) \cup \left(\bigcup_{a < \delta} (M\mathbb{H}_1^a) + \varepsilon \mathbb{B}_1\right).$$

By Lemma 4.7 the set $B$ contains the set $M\mathbb{H}_1^a + \varepsilon\mathbb{B}_1$ for any $a \in [\delta, r]$. This is true also for $a < \delta$, trivially, by the definition of $B$. Consequently, by Borell's inequality (see [5] or Theorem 5.1 in [48]), for any $a \leq r$,

$$P(W^a \notin B) \leq P(W^a \notin M\mathbb{H}_1^a + \varepsilon\mathbb{B}_1) \leq 1 - \Phi(\Phi^{-1}(e^{-\phi_0^a(\varepsilon)}) + M)$$
$$\leq 1 - \Phi(\Phi^{-1}(e^{-\phi_0^r(\varepsilon)}) + M),$$

because $e^{-\phi_0^a(\varepsilon)} = P(\sup_{t \in aT} |W_t| \leq \varepsilon)$ is decreasing in $a$. For

$$M \geq -2\Phi^{-1}(e^{-\phi_0^r(\varepsilon)}),$$

the right-hand side is bounded by $1 - \Phi(M/2) \leq e^{-M^2/8}$. The latter condition is certainly satisfied if (cf. Lemma 4.10),

$$M \geq 4\sqrt{\phi_0^r(\varepsilon)} \quad \text{and} \quad e^{-\phi_0^r(\varepsilon)} < 1/4.$$

Here $e^{-\phi_0^r(\varepsilon)} \leq e^{-\phi_0^1(\varepsilon)}$ for $r > 1$ and is certainly smaller than $1/4$ if $\varepsilon$ is smaller than some fixed $\varepsilon_1$. Therefore, in view of Lemma 4.6 the inequalities are satisfied if

$$(5.2) \qquad M^2 \geq 16 C_4 r^d (\log(r/\varepsilon))^{1+d}, \qquad r > 1, \qquad \varepsilon < \varepsilon_1 \wedge \varepsilon_0.$$

In view of Lemma 4.9, for $r$ larger than a positive constant depending on $d$ and the density of $A$ only,

$$(5.3) \qquad P(W^A \notin B) \leq P(A > r) + \int_0^r P(W^a \notin B) g(a)\, da$$
$$\leq \frac{2C_2 r^{p-d+1} e^{-D_2 r^d \log^q r}}{D_2 d \log^q r} + e^{-M^2/8}.$$

This inequality is true for any $B = B_{M,r,\delta,\varepsilon}$ with $M, r, \delta, \varepsilon$ satisfying (5.2).

By Lemma 4.5, for $M\sqrt{r/\delta} > 2\varepsilon$ and $r > a_0$,

$$\log N\left(2\varepsilon, M\sqrt{\frac{r}{\delta}}\mathbb{H}_1^r + \varepsilon\mathbb{B}_1, \|\cdot\|_\infty\right) \leq \log N\left(\varepsilon, M\sqrt{\frac{r}{\delta}}\mathbb{H}_1^r, \|\cdot\|_\infty\right)$$
$$\leq K r^d \left(\log\left(\frac{M\sqrt{r/\delta}}{\varepsilon}\right)\right)^{1+d}.$$

By Lemma 4.8 every element of $M\mathbb{H}_1^a$ for $a < \delta$ is within uniform distance $\delta\sqrt{d}\tau M$ of a constant function for a constant in the interval $[-E, E]$, for $E = M\sqrt{\|\mu\|}$. It follows that, for $\varepsilon > \delta\sqrt{d}\tau M$,

$$N\left(3\varepsilon, \bigcup_{a<\delta}(M\mathbb{H}_1^a) + \varepsilon\mathbb{B}_1, \|\cdot\|_\infty\right) \leq N(\varepsilon, [-E, E], |\cdot|) \leq \frac{2E}{\varepsilon}.$$

The covering number of a union is bounded by the sum of the covering numbers. Therefore, with the choice $\delta = \varepsilon/(2\sqrt{d}\tau M)$, together the last two displays yield, since $\log(x+y) \leq \log(2(x \vee y)) \log x + 2 \log y$ for $x \geq 1, y \geq 2$, for $2E/\varepsilon \geq 2$,

$$\log N(3\varepsilon, B, \|\cdot\|_\infty) \leq K r^d \left(\log\left(\frac{M^{3/2}\sqrt{2\tau r}d^{1/4}}{\varepsilon^{3/2}}\right)\right)^{1+d}$$

(5.4)
$$+ 2 \log \frac{2M\sqrt{\|\mu\|}}{\varepsilon}.$$

This inequality is valid for any $B = B_{M,r,\delta,\varepsilon}$ with $\delta = \varepsilon/(2\sqrt{d}\tau M)$, and any $M, r, \varepsilon$ with

(5.5) $\quad M^{3/2}\sqrt{2\tau r}d^{1/4} > 2\varepsilon^{3/2}, \qquad r > a_0, \qquad M\sqrt{\|\mu\|} > \varepsilon.$

In the remainder of the proof we make special choices for these parameters, depending on the assumption on $w_0$.

5.1. *Hölder smoothness.* Suppose that $w_0 \in C^\alpha[0,1]^d$ for some $\alpha > 0$. In view of Lemmas 4.3 and 4.6, for every $a_0$ there exist positive constants $\varepsilon_0 < 1/2$, $C$, $D$ and $K$ that depend on $w$ and $\mu$ only such that, for $a > a_0$, $\varepsilon < \varepsilon_0$ and $\varepsilon > Ca^{-\alpha}$,

$$\phi^a_{w_0}(\varepsilon) \leq Da^d + C_4 a^d \left(\log \frac{a}{\varepsilon}\right)^{1+d} \leq K_1 a^d \left(\log \frac{a}{\varepsilon}\right)^{1+d}$$

for $K_1$ depending on $a_0$, $\mu$ and $d$ only. Therefore, for $\varepsilon < \varepsilon_0 \wedge Ca_0^{-\alpha}$ [so that $(C/\varepsilon)^{1/\alpha} > a_0$], by (5.1),

$$P(\|W^A - w_0\|_\infty \leq 2\varepsilon) \geq \int_0^\infty e^{-\phi^a_{w_0}(\varepsilon)} g(a)\, da$$

$$\geq \int_{(C/\varepsilon)^{1/\alpha}}^{2(C/\varepsilon)^{1/\alpha}} e^{-K_1 a^d \log^{1+d}(a/\varepsilon)} g(a)\, da$$

$$\geq C_1 e^{-K_2(1/\varepsilon)^{d/\alpha}(\log(1/\varepsilon))^{(1+d)\vee q}} \left(\frac{C}{\varepsilon}\right)^{p/\alpha} \left(\frac{C}{\varepsilon}\right)^{1/\alpha},$$

in view of (3.4), for a constant $K_2$ that depends only on $K_1, C, D_1, d, \alpha, q$. We conclude that $P(\|W^A - w_0\|_\infty \leq \varepsilon_n) \geq \exp(-n\varepsilon_n^2)$ for $\varepsilon_n$ a large multiple of $n^{-1/(2+d/\alpha)}(\log n)^\gamma$, for $\gamma = ((1+d) \vee q)/(2+d/\alpha)$, and sufficiently large $n$.

By (5.2)–(5.3) $P(W^A \notin B)$ is bounded above by a multiple of $\exp(-C_0 n\varepsilon_n^2)$ for an arbitrarily large constant $C_0$ if (5.2) holds and

(5.6)
$$D_2 r^d (\log r)^q \geq 2C_0 n\varepsilon_n^2,$$
$$r^{p-d+1} \leq e^{C_0 n\varepsilon_n^2},$$
$$M^2 \geq 8C_0 n\varepsilon_n^2.$$

Given $C_0$ we first choose $r = r_n$ as the minimal solution to the first equation, and next we choose $M = M_n$ to satisfy the third equation and (5.2). The second equation is then automatically satisfied, for large $n$.

With these choices of $M$ and $r$ and $\bar{\varepsilon}_n$ bounded below by a power of $n$ the right-hand side of (5.4) is bounded by a multiple of $r_n^d (\log n)^{1+d} + \log n$. This is bounded by $n\bar{\varepsilon}_n^2$ for $\bar{\varepsilon}_n^2$ a large multiple of $(r_n^d/n)(\log n)^{1+d}$. Inequalities (5.5) are clearly satisfied.

5.2. *Infinite smoothness, $r \geq \nu$.* Suppose that $w_0$ is the restriction of a function $w_0 \in \mathcal{A}^{\gamma,r}(\mathbb{R}^d)$ for $r \geq \nu$, and that the spectral density is bounded below by a multiple of $\exp(-D_3 \|\lambda\|^\nu)$ for some positive constants $D_3$ and $\nu$. By combining the first part of Lemma 4.4 and Lemma 4.6, we see that there exist positive constants $a_0 < a_1$, $\varepsilon_0$, $K_1$ and $C_4$ that depend on $w$ and $\mu$ only such that, for $a \in [a_0, a_1]$ and $\varepsilon < \varepsilon_0$,

$$\phi_{w_0}^a(\varepsilon) \leq K_1 + C_4 a^d \left(\log \frac{a}{\varepsilon}\right)^{1+d}.$$

Consequently, by (5.1),

$$P(\|W^A - w_0\|_\infty \leq 2\varepsilon) \geq \int_0^\infty e^{-\phi_{w_0}^a(\varepsilon)} g(a)\, da$$

$$\geq e^{-K_1 - C_4 a_1^d \log^{1+d}(a_1/\varepsilon)} P(a_0 < A < a_1).$$

We conclude that $P(\|W^A - w_0\|_\infty \leq \varepsilon_n) \geq \exp(-n\varepsilon_n^2)$ for $\varepsilon_n$ a large multiple of $n^{-1/2}(\log n)^{(d+1)/2}$, and sufficiently large $n$.

Next we choose $B$ of the form as before, with $r$ and $M$ solving (5.6) and satisfying (5.2), that is, $r_n^d$ and $M_n^2$ large multiples of $(\log n)^{d+1}$. Then (5.2)–(5.3) show that $P(W^A \notin B)$ is bounded above by a multiple of $\exp(-C_0 n\varepsilon_n^2)$, and the right-hand side of (5.4) is bounded by a multiple of $r_n^d(\log(1/\varepsilon) + \log \log n)^{1+d} + \log(1/e) + \log \log n$. For $\varepsilon = \bar{\varepsilon}_n$ a large multiple of $n^{-1/2}(\log n)^{d+1}$ this is bounded above by $n\bar{\varepsilon}_n^2$.

5.3. *Infinite smoothness, $r < \nu$.* Consider the situation as in the preceding section, but now with $r < \nu$. Combining the second part of Lemma 4.4 and Lemma 4.6, we see that there exist positive constants $a_0$, $\varepsilon_0$, $C$, $D$, $K_1$ and $C_4$ that depend on $w$ and $\mu$ only and $\gamma' > \gamma$ such that, for $a > a_0$, $\varepsilon < \varepsilon_0$ and $C \exp(-\gamma' a^r) < \varepsilon$,

$$\phi_{w_0}^a(\varepsilon) \leq D a^d + C + 4 a^d \left(\log \frac{a}{\varepsilon}\right)^{1+d}.$$

Consequently, by (5.1), for constants $D_1, D_2$ that depend on $w$ and $\mu$ only,

$$P(\|W^A - w_0\|_\infty \leq 2\varepsilon) \geq \int_{(\log(C/\varepsilon)/\gamma')^{1/r}}^\infty e^{-\phi_{w_0}^a(\varepsilon)} g(a)\, da$$

$$\geq D_2 e^{-D_1(\log(1/\varepsilon))^{d/r+d+1}}.$$

We conclude that $P(\|W^A - w_0\|_\infty \leq \varepsilon_n) \geq \exp(-n\varepsilon_n^2)$ for $\varepsilon_n$ a large multiple of $n^{-1/2}(\log n)^{d/(2r)+(d+1)/2}$, and sufficiently large $n$.

Next we choose $B$ of the form as before, with $r$ and $M$ solving (5.6), that is, $r_n^d$ and $M_n^2$ large multiples of $(\log n)^{d/r+d+1}$. Then (5.2) and (5.3) show that $P(W^A \notin B)$ is bounded above by a multiple of $\exp(-C_0 n\varepsilon_n^2)$, and the right-hand side of (5.4) is bounded by a multiple of $r_n^d(\log(1/\varepsilon) + \log\log n)^{1+d} + \log(1/\varepsilon) + \log\log n$. For $\varepsilon = \bar\varepsilon_n$ a large multiple of $n^{-1/2}(\log n)^{d+1+d/(2r)}$ this is bounded above by $n\bar\varepsilon_n^2$.

DEPARTMENT OF MATHEMATICS
VRIJE UNIVERSITEIT
DE BOELELAAN 1081A
1081 HV AMSTERDAM
THE NETHERLANDS
E-MAIL: aad@cs.vu.nl
　　　　harry@cs.vu.nl